\documentclass[reqno]{amsart}
\usepackage{amsthm}
\usepackage{amsmath}
\usepackage{enumerate}
\usepackage{graphicx}
\usepackage{alphabeta}

\newcommand{\fake}[0]{\fontfamily{ptm}\digamma}
\newcommand{\Fake}[0]{F}
\newcommand{\abs}[1]{\left\lvert#1\right\rvert}
\newcommand{\floor}[1]{\left\lfloor#1\right\rfloor}
\DeclareMathOperator{\Res}{Res}
\newtheorem{thm}{Theorem}
\newcommand{\figdir}[0]{Figs}
\newcommand{\gfextn}[0]{jpeg}

\begin{document}

\title{Fake Mu's}

\author[G.~Martin]{Greg Martin}
\address{Department of Mathematics\\ University of British Columbia\\
Vancouver, BC, Canada}
\email{gerg@math.ubc.ca}

\author[M.~J. Mossinghoff]{Michael J. Mossinghoff}
\address{Center for Communications Research\\
Princeton, NJ, USA}
\email{m.mossinghoff@idaccr.org}

\author[T.~S. Trudgian]{Timothy S. Trudgian}\thanks{This work was supported by a Future Fellowship (FT160100094 to T.~S. Trudgian) from the Australian Research Council.}
\address{School of Science\\
UNSW Canberra at ADFA\\
ACT, Australia}
\email{t.trudgian@adfa.edu.au}

\date{\today}

\begin{abstract}
Let $\fake(n)$ denote a multiplicative function with range $\{-1,0,1\}$, and let $\Fake(x) = \sum_{n=1}^{\floor{x}} \fake(n)$.
Then $\Fake(x)/\sqrt{x} = a\sqrt{x} + b + E(x)$, where $a$ and $b$ are constants and $E(x)$ is an error term that either tends to $0$ in the limit, or is expected to oscillate about $0$ in a roughly balanced manner.
We say $\Fake(x)$ has \textit{persistent bias} $b$ (at the scale of $\sqrt{x}$) in the first case, and \textit{apparent bias} $b$ in the latter.
For example, if $\fake(n)=\mu(n)$, the M\"{o}bius function, then $\Fake(x) = \sum_{n=1}^{\floor{x}} \mu(n)$ has $b=0$ so exhibits no persistent or apparent bias, while if $\fake(n)=\lambda(n)$, the Liouville function, then $\Fake(x) = \sum_{n=1}^{\floor{x}} \lambda(n)$ has apparent bias $b=1/\zeta(1/2)$.
We study the bias when $\fake(p^k)$ is independent of the prime $p$, and call such functions \textit{fake $\mu's$}.
We investigate the conditions required for such a function to exhibit a persistent or apparent bias, determine the functions in this family with maximal and minimal bias of each type, and characterize the functions with no bias of either type.
For such a function $\Fake(x)$ with apparent bias $b$, we also show that $\Fake(x)/\sqrt{x}-a\sqrt{x}-b$ changes sign infinitely often.
\end{abstract}

\maketitle

\section{Introduction}

Let $M(x)$ denote the \textit{Mertens function}, defined as the sum of the values of the M\"obius function $\mu(n)$ over all positive integers $n\leq x$,
\[
M(x) = \sum_{n=1}^{\floor{x}} \mu(n).
\]
This function has a long history in number theory.
In 1897 Mertens \cite{Mertens} tabulated its values up to $10^4$ and ventured that it was ``highly probable'' that $\abs{M(x)}\leq\sqrt{x}$ for all $x\geq1$.
It is well known that a bound of this form would imply the Riemann hypothesis (RH) and the simplicity of the zeros of the Riemann zeta function (SZ).
In fact, both of those statements would follow if $M(x)/\sqrt{x}$ were bounded either above or below by any fixed constant.
In an influential paper of 1942, Ingham \cite{Ingham1942} proved that any such bound on $M(x)/\sqrt{x}$ would imply significantly more---it would also follow that there exist infinitely many integer relations among the ordinates of the zeros of the zeta function in the upper half plane.
Mertens' conjecture remained open until Odlyzko and te Riele disproved it in 1985 \cite{OTR}; see \cite{Hurst,Ng} for more information on oscillations of $M(x)/\sqrt{x}$.

In a similar manner, let $L(x)$ denote the sum of the values of the Liouville function $\lambda(n)$ over positive integers $n\leq x$,
\[
L(x) = \sum_{n=1}^{\floor{x}} \lambda(n).
\]
Recall that $\lambda(n)=(-1)^{\Omega(n)}$, where $\Omega(n)$ counts the number of prime divisors of $n$, including multiplicities.
In 1919 P\'olya \cite{Polya} computed values of $L(x)$ up to about $x=1500$, and noted that $L(x)$ was never positive over this range, once $x\geq2$.
P\'olya remarked that if $L(x)$ never changed sign for sufficiently large $x$, then the Riemann hypothesis would follow, as well as the simplicity of the zeros of the zeta function.
Ingham's work also covered the function $L(x)$, so it was known since then as well that bounding $L(x)/\sqrt{x}$ either above or below would also imply the existence of integer relations among the nontrivial zeros of the zeta function.
Haselgrove proved that $L(x)$ does change sign infinitely often in 1958 \cite{Haselgrove}, and specific values for sign changes were found later \cite{Borwein,Lehman,TanakaLiouville}.
Additional information concerning oscillations in $L(x)/\sqrt{x}$ can be found in \cite{HumphriesJNT,MT2017}.

\begin{figure}[tb]
\begin{center}
\includegraphics[width=4in]{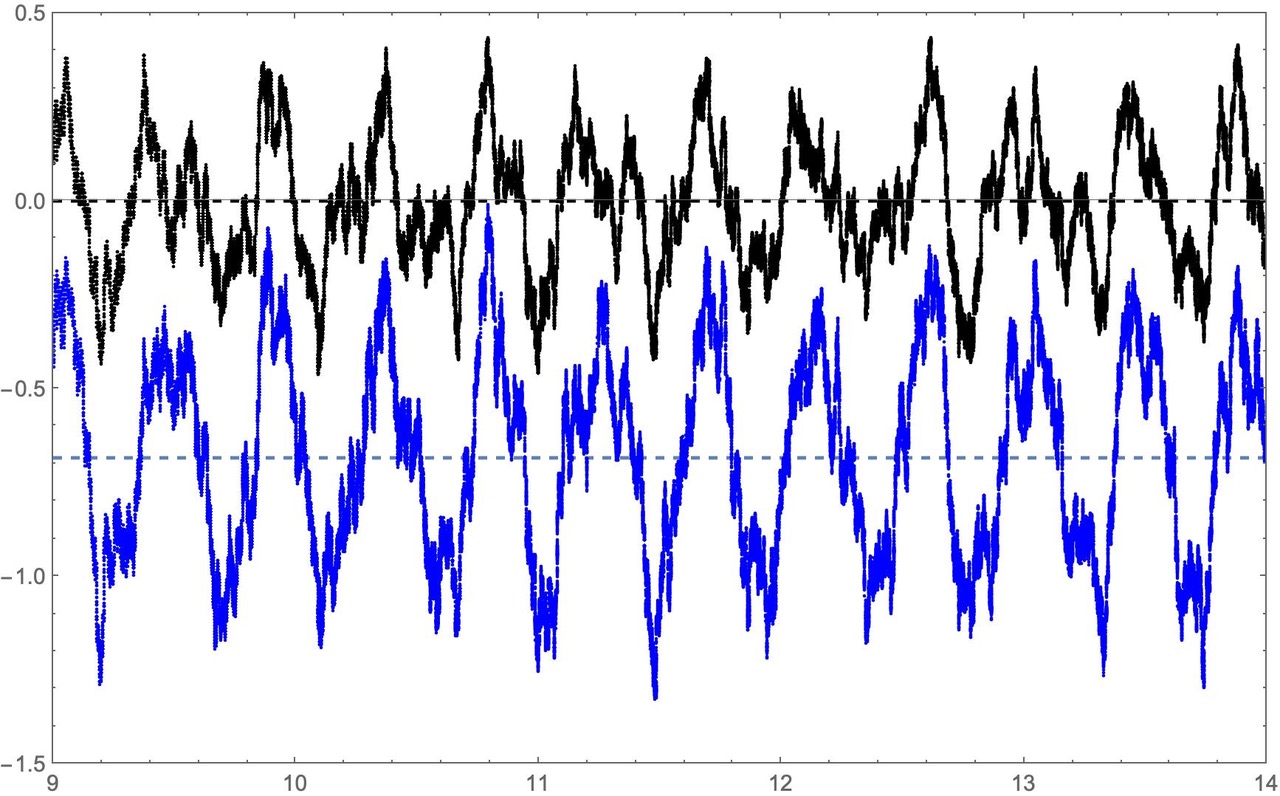}
\end{center}
\caption{Normalized Mertens function $e^{-u/2}M(e^u)$, centered at $0$ (top), and P\'olya function $e^{-u/2}L(e^u)$, centered at $1/\zeta(1/2)=-0.68476\ldots$ (bottom).}\label{figML}
\end{figure}

It is interesting that $M(x)$ changes sign rather frequently, while $L(x)$ seems skewed toward negative values.
This is evident in the plots of $M(x)/\sqrt{x}$ and $L(x)/\sqrt{x}$ shown in Figure~\ref{figML}.
A reason for this contrast can be seen by analyzing the Dirichlet series for $\mu(n)$ and $\lambda(n)$,
\[
\sum_{n\geq1} \frac{\mu(n)}{n^s} = \frac{1}{\zeta(s)}
\textrm{\quad and\quad}
\sum_{n\geq1} \frac{\lambda(n)}{n^s} = \frac{\zeta(2s)}{\zeta(s)}.
\]
Assuming RH and SZ, by applying Perron's formula \cite[Lem.\ 3.12]{Touchy} to these Dirichlet series one obtains
\begin{align*}
\frac{M(x)}{\sqrt{x}} &= \sum_{\substack{\rho_n=1/2+i\gamma_n\\\abs{\gamma_n}\leq T}} \frac{x^{i\gamma_n}}{\rho_n\zeta'(\rho_n)} + E_{1}(x,T),\\
\frac{L(x)}{\sqrt{x}} &= \frac{1}{\zeta(1/2)} + \sum_{\substack{\rho_n=1/2+i\gamma_n\\\abs{\gamma_n}\leq T}} \frac{\zeta(2\rho_n) x^{i\gamma_n}}{\rho_n\zeta'(\rho_n)} + E_{2}(x,T),
\end{align*}
where $E_{1}(x, T)$ and $E_{2}(x,T)$ denote  error terms that tend to $0$ as $T\to\infty$.
If we model the $x^{i\gamma_n}$ as random unit vectors in the complex plane, then we would expect $L(x)/\sqrt{x}$ to oscillate about the constant term $1/\zeta(1/2)=-0.68476523\ldots$\,, and $M(x)/\sqrt{x}$ to oscillate about $0$, since there is no constant term in that case.
The constant term $1/\zeta(1/2)$ in the expression for $L(x)$ arises from the pole in the Dirichlet series for the Liouville function lying at $s=1/2$; the series for the M\"obius function exhibits no such pole.

\begin{figure}[tb]
\begin{center}
\includegraphics[width=4in]{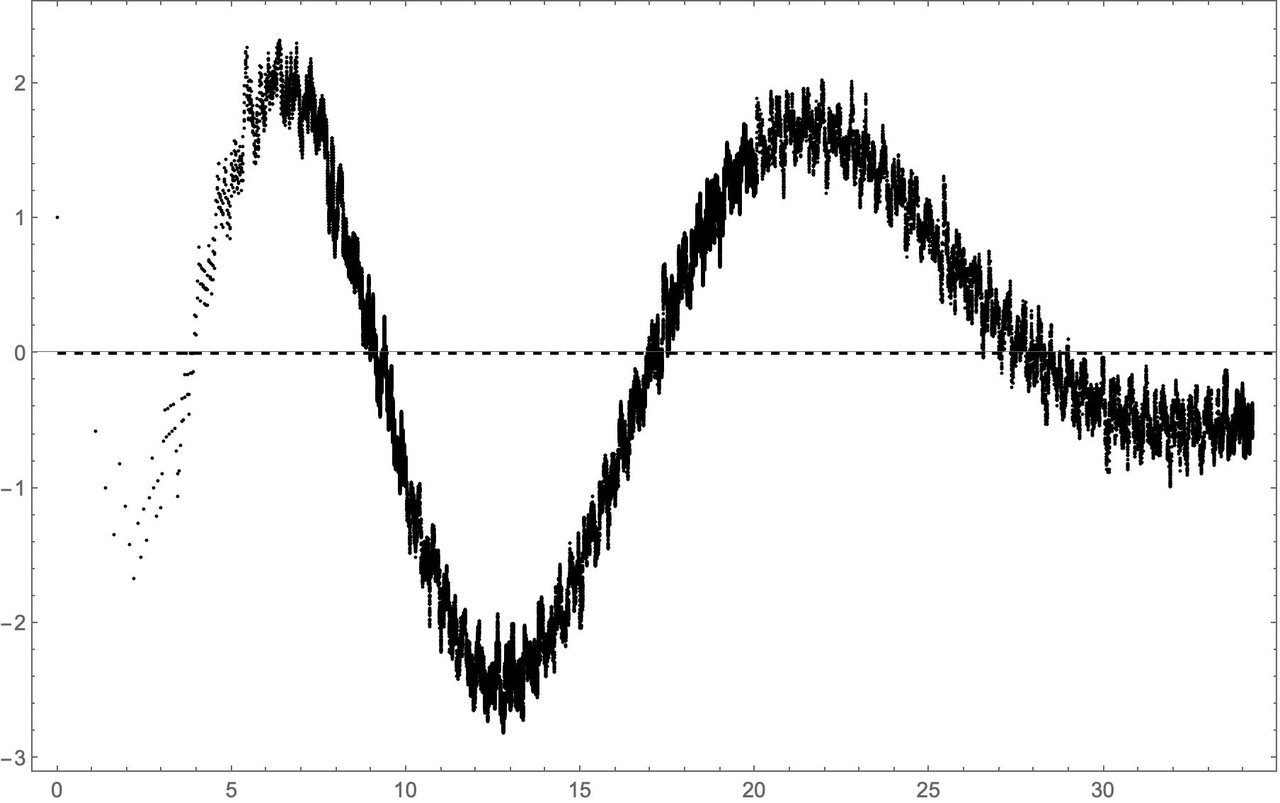}
\end{center}
\caption{The normalized function $e^{-u/2}\Xi(e^u)$ for $u\leq\log(1.5\cdot10^{15})$.}\label{figH1Q}
\end{figure}

One may also consider the similar function $\xi(n)=(-1)^{\omega(n)}$, where $\omega(n)$ denotes the number of distinct prime factors of $n$.
Let $\Xi(x)$ denote its summatory function,
\[
\Xi(x) = \sum_{n=1}^{\floor{x}} \xi(n).
\]
This function was recently investigated by the second and third authors \cite{Dickens,MTIJNT}, who showed that the corresponding Dirichlet series exhibits no pole at $s=1/2$, so one expects $\Xi(x)/\sqrt{x}$ to be unbiased in sign for large $x$.
(This was also shown implicitly by van de Lune and Dressler \cite{Dressler}.)
A plot of $\Xi(x)/\sqrt{x}$ on a logarithmic scale appears in Figure~\ref{figH1Q}.
The macroscopic shape here arises from the contributions of poles of the corresponding Dirichlet series on vertical lines with real part $\sigma\leq1/4$ (assuming RH and SZ), and these contributions dampen as $x$ grows large.
Here and throughout this article, we write $s=\sigma+i t$ with $\sigma$ and $t$ real.

In this article we consider some families of arithmetic functions similar to $\mu(n)$, $\lambda(n)$, and $\xi(n)$, and measure an expected bias in their summatory functions by computing the poles in their Dirichlet series.
The functions $\fake(n)$ we consider are defined by setting $\fake(1)=1$ and
\[
\fake(p^k) = \epsilon_k \in \{-1,0,1\}
\]
for each $k\geq1$ and each prime $p$, where $\epsilon_k$ does not depend on $p$, and requiring that $\fake(n)$ be multiplicative.
This family clearly includes $\mu(n)$, $\lambda(n)$, and $\xi(n)$, as well as other functions commonly studied in number theory, such as $\mu^2(n)$, the indicator function for squarefree integers.
We refer to this family of functions $\fake(n)$ as \textit{fake $\mu$'s}, since they generalize the M\"obius function in particular.
(Here, $\fake$ is the archaic Greek character digamma.)

For such a function $\fake(n)$, we denote its summatory function by
\[
\Fake(x) = \sum_{n=1}^{\floor{x}} \fake(n),
\]
and we wish to investigate whether this function is expected to exhibit a bias at a particular scale.
Unlike $M(x)$, $L(x)$, or $\Xi(x)$, some of the functions in the family we study have a linear term that we must discount before normalizing.
For example, it is well known that $\sum_{n=1}^{\floor{x}} \mu^2(n) \sim x/\zeta(2)$, and if we set $\fake(p^k)=1$ for all $k\geq0$ then the resulting summatory function is simply $\sum_{n=1}^{\floor{x}} 1 = \floor{x}$.
For a particular choice of $\fake(n)$, we let $a$ denote the limiting value of $\Fake(x)/x$ as $x\to\infty$, so this is the residue of the Dirichlet series $\sum_{n\geq1} \fake(n)/n^s$ at $s=1$ when a pole occurs here.
Then we are interested in the behavior of the normalized summatory function
\begin{equation}\label{eqnNormalize}
\frac{\Fake(x)-ax}{\sqrt{x}}.
\end{equation}
If this expression has a nonzero limit $b$ as $x\to\infty$, then we say $\Fake(x)$ has \textit{persistent bias} $b$.
If it is expected to oscillate about a value $b$ in a roughly balanced way for large $x$, then we say $b$ is the \textit{apparent bias} of $\Fake(x)$.
For example, $L(x)$ has apparent bias $1/\zeta(1/2)$, while $M(x)$ and $\Xi(x)$ have no persistent or apparent bias.

We remark that when we say a function has no persistent or apparent bias, we mean this only at the scale of $\sqrt{x}$.
A function with no bias at this scale could well see one at a smaller scale, for example, if its Dirichlet series has a pole on the real axis at some $\sigma<1/2$.
Similarly, a function could exhibit oscillations at a smaller scale, if its Dirichlet series has poles at complex values $\sigma+it$ with $\sigma<1/2$.
This occurs for example when $\fake(n)=\mu^2(n)$, which has poles on the line $\sigma=1/4$, assuming RH\@.
(See \cite{MOST} for more on oscillations in its summatory function.)

One natural subset of the family of functions we consider here was studied by Tanaka in 1980 \cite{Tanaka}.
For a positive integer $r$ let
\[
\mu_r(p^k) = \begin{cases}
(-1)^k, & k<r,\\
0, & k\geq r,
\end{cases}
\]
so that $\mu_2(n)=\mu(n)$ and $\mu_\infty(n)=\lambda(n)$.
The Dirichlet series for $\mu_r(n)$ with $r\geq1$ is given by
\[
T_{\mu_r}(s) = \begin{cases}
\displaystyle\frac{\zeta(2s)\zeta(rs)}{\zeta(s)\zeta(2rs)}, & \textrm{$r$ odd},\\[2ex]
\displaystyle\frac{\zeta(2s)}{\zeta(s)\zeta(rs)}, & \textrm{$r$ even}.\\
\end{cases}
\]
Each of these functions is analytic in the half-plane $\sigma>1/2$ assuming RH, and the residue of $T_{\mu_r}(s)$ at $s=1/2$ is
\[
\Res_{1/2}(T_{\mu_r}) =  \begin{cases}
\displaystyle\frac{1}{2\zeta(1/2)\zeta(r/2)}, & \textrm{$r\geq3$ odd},\\[2ex]
\displaystyle\frac{\zeta(r/2)}{2\zeta(1/2)\zeta(r)}, & \textrm{$r\geq4$ even},
\end{cases}
\]
so the apparent bias of the summatory function $M_r(x)$ of $\mu_r(n)$ is $b_r = 2\Res_{1/2}(T_{\mu_r})$.
Tanaka showed that for each $r\geq2$ the function $M_r(x)/\sqrt{x}-b_r$ changes sign infinitely often, proving in fact that
\begin{equation}\label{eqnTanaka}
\liminf_{x\to\infty} \frac{M_r(x)}{\sqrt{x}} - b_r < 0
\quad\textrm{and}\quad
\limsup_{x\to\infty} \frac{M_r(x)}{\sqrt{x}} - b_r > 0.
\end{equation}
We note that each of the values $b_r$ with $r\geq3$ is negative, and the largest in absolute value among them occurs for $r=3$, where the apparent bias is more than twice that of $L(x)$:
\[
b_3 = \frac{1}{\zeta(1/2)\zeta(3/2)} = -1.488169\ldots\,.
\]
This provides an initial benchmark for our work: do there exist other functions in our family with an even more pronounced negative bias?
More generally, we seek to answer the following questions.
\begin{itemize}
\item What features of $\fake(n)$ allow us to move from no apparent bias in the summatory function, as with $M(x)$, to an apparent bias, as with $L(x)$?
\item The apparent biases in $L(x)$ and $M_r(x)$ with $r\geq3$ are negative.
Do there exist arithmetic functions in this family with positive apparent bias?
\item How large of a bias exists for these functions, whether apparent or persistent, in the positive or negative directions, and what functions $\fake(n)$ attain the extremal values?
\item Can one characterize the unbiased functions in the families we consider?
\end{itemize}

We answer these questions in this article.
We first consider the family of functions $\fake(n)$ with range $\{-1,1\}$, so where $\fake(p^k)=\pm1$ for each $k$.
This allows interpolating $\lambda(n)$ and $\xi(n)$ in particular, and studying how the apparent bias behaves over these functions.
For this case, we prove the following theorem in Section~\ref{secpm1}.

\begin{thm}\label{thmpm1}
Suppose $\fake(n)$ is a multiplicative function satisfying $\fake(p^k)=\epsilon_k\in\{-1,1\}$ for each prime $p$, and let $\Fake(x)$ denote its summatory function.
Then
\begin{enumerate}[(i)]
\item If $\epsilon_1=1$, or if $\epsilon_1=\epsilon_2=-1$, then $\Fake(x)$ has no persistent or apparent bias.
\item\label{item1b} If $\epsilon_1=-1$ and $\epsilon_2=1$, then $\Fake(x)$ has apparent bias $b\in[A_1,B_1]$, where
\begin{align*}
A_1 &= \frac{1}{\zeta(1/2)}\prod_p\left(1+\frac{2}{p(\sqrt{p}-1)}\right) = -10.29174\ldots,\\
B_1 &= \frac{1}{\zeta(1/2)}\prod_p\left(1-\frac{2}{p^{3/2}(\sqrt{p}-1)}\right) = 0.16918\ldots,
\end{align*}
and both extreme values are achieved.
\item\label{item1c} There exist uncountably many functions with no persistent or apparent bias having $\epsilon_1=-1$ and $\epsilon_2=1$, and these functions may be characterized precisely.
\end{enumerate}
\end{thm}

We then consider the broader case where $\fake(p^k)\in\{-1,0,1\}$ and prove the following theorem in Section~\ref{seczpm1}.

\begin{thm}\label{thmzpm1}
Suppose $\fake(n)$ is a multiplicative function satisfying $\fake(p^k)=\epsilon_k\in\{-1,0,1\}$ for each prime $p$, and let $\Fake(x)$ denote its summatory function.
Then
\begin{enumerate}[(i)]
\item If $\epsilon_1=1$, or $\epsilon_1,\epsilon_2\in\{0,-1\}$, then $\Fake(x)$ has no persistent or apparent bias.
\item\label{item2b} If $\epsilon_1=-1$ and $\epsilon_2=1$, then $\Fake(x)$ has apparent bias $b\in[A_1,B_1]$ from Theorem~\ref{thmpm1}.
Further, the functions with no apparent bias in this family may be characterized precisely.
\item If $\epsilon_1=0$ and $\epsilon_2=1$, then $\Fake(x)$ has persistent bias $b\in[A_2,B_2]$, where
\begin{align*}
A_2 &= \prod_p\left(1-\frac{1}{p^{3/2}}-\frac{2}{p^2}\right) = 0.051526\ldots,\\
B_2 &= \frac{\zeta(3)}{\zeta(3/2)} = 2.17325\ldots,
\end{align*}
and both extreme values are achieved.
\end{enumerate}
\end{thm}

We also prove that the normalized summatory functions $\Fake(x)/\sqrt{x}$ for the fake $\mu$'s from part~(\ref{item2b}) of Theorems~\ref{thmpm1} and~\ref{thmzpm1} are infinitely often larger than their apparent bias $b$ as $x\to\infty$, and infinitely often smaller than this value as $x\to\infty$.
We prove the following result in Section~\ref{secBias} by generalizing the argument of Tanaka that produced \eqref{eqnTanaka} for the functions $\mu_r(n)$.

\begin{thm}\label{thmBias}
Suppose $\fake(n)$ is a multiplicative function with $\fake(p)=-1$, $\fake(p^2)=1$, and $\fake(p^k)=\epsilon_k\in\{-1,0,1\}$ for $k\geq3$, for each prime $p$.
Let $\Fake(x)$ denote its summatory function, and let $b\in[A_1,B_1]$ from Theorem~\ref{thmpm1} denote its apparent bias.
Then
\[
\liminf_{x\to\infty} \frac{\Fake(x)}{\sqrt{x}} < b
\quad\textrm{and}\quad
\limsup_{x\to\infty} \frac{\Fake(x)}{\sqrt{x}} > b.
\]
\end{thm}

\begin{figure}[tb]
\begin{center}
\includegraphics[width=4in]{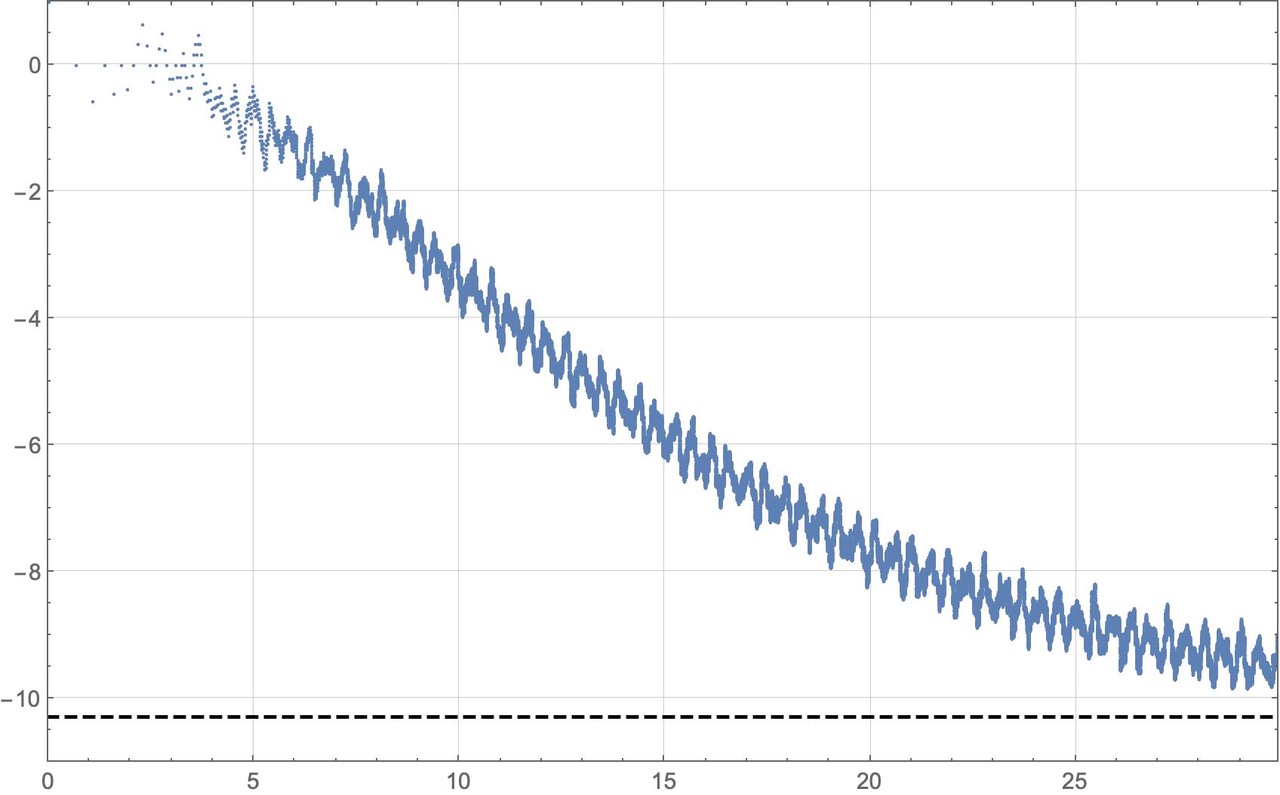}
\end{center}
\caption{Normalized summatory function $e^{-u/2}\Fake_{\textrm{min}}(e^u)$ with $u\leq\log(10^{13})$, for $\fake_{\textrm{min}}(n)$ from \eqref{eqnfakemin}, with minimal apparent bias.}\label{figMaxNegBias}
\end{figure}

\begin{figure}[tb]
\begin{center}
\includegraphics[width=4in]{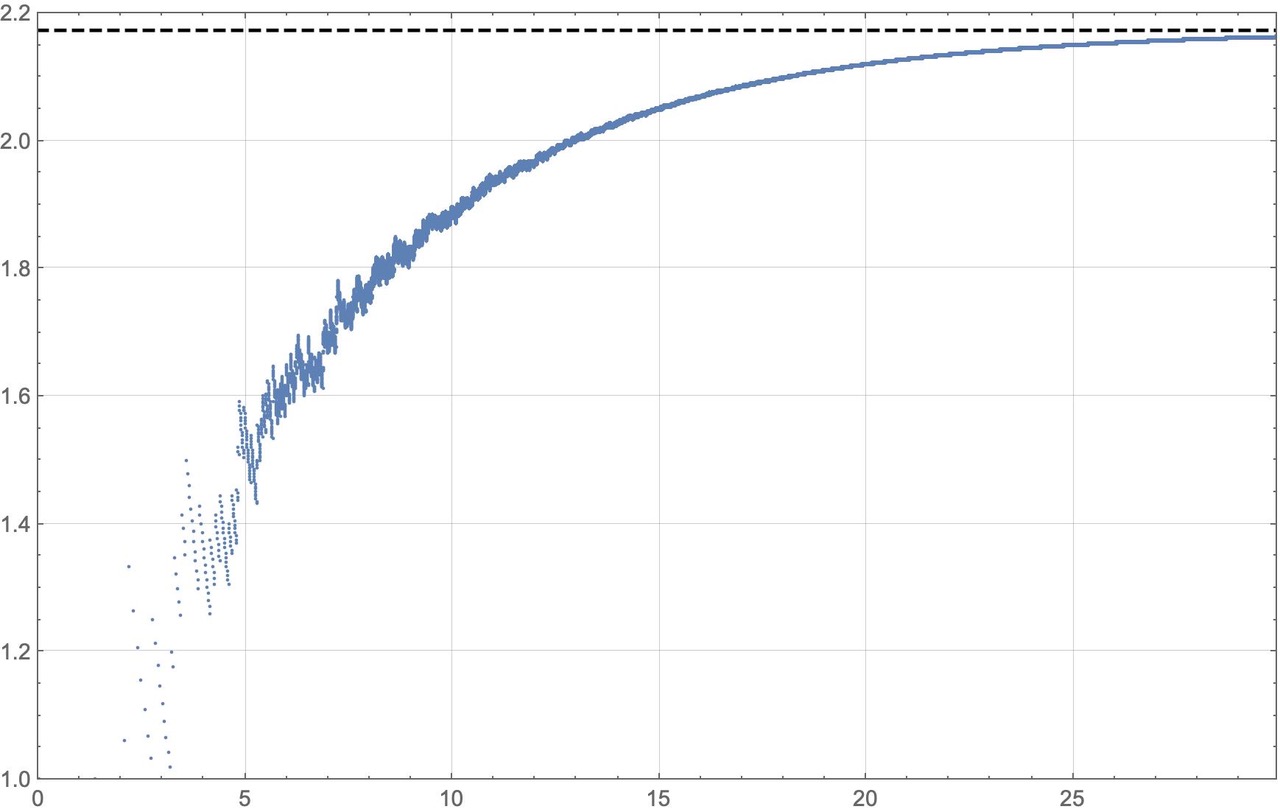}
\end{center}
\caption{Normalized summatory function $e^{-u/2}\Fake_{\textrm{Max}}(e^u)$ with $u\leq\log(10^{13})$, for $\fake_{\textrm{Max}}(n)$ from \eqref{eqnfakeMax}, with maximal persistent bias.}\label{figMaxPosBias}
\end{figure}

The most negatively biased function among the fake $\mu$'s considered here, where $\Fake(x)/\sqrt{x}$ has apparent bias $A_1=-10.29\ldots$ from Theorem~\ref{thmpm1}, is given by evaluating the M\"obius function on the power-free part of the argument:
\begin{equation}\label{eqnfakemin}
\fake_{\textrm{min}}(p^k) = \begin{cases}
-1, & k=1,\\
1, & \textrm{otherwise}.
\end{cases}
\end{equation}
The most positively biased function in our study, where $\Fake(x)/\sqrt{x}$ has persistent bias $B_2=2.17\ldots$ from Theorem~\ref{thmzpm1},  arises from a function that indicates whether its argument $n$ is \textit{powerful}, that is, when every prime divisor $p$ of $n$ satisfies $p^2\mid n$:
\begin{equation}\label{eqnfakeMax}
\fake_{\textrm{Max}}(p^k) = \begin{cases}
0, & k =1,\\
1, & \textrm{otherwise}.
\end{cases}
\end{equation}
Plots of the normalized summatory functions for these two extremal examples are shown in Figures~\ref{figMaxNegBias} and~\ref{figMaxPosBias} for $x\leq10^{13}$, displayed on a logarithmic scale.
Plots for other functions relevant to Theorems~\ref{thmpm1} and~\ref{thmzpm1} are also displayed later in the article.

We remark that biases and asymptotic behavior of other generalizations of $\mu(n)$ and $\lambda(n)$ occur in the literature.
For example, Humphries, Shekatkar, and Woo \cite{HSW} recently studied the summatory function of $\lambda(n; q,a)=(-1)^{\Omega(n;q,a)}$, where $\Omega(n;q,a)$ denotes the number of prime factors $p$ of $n$ (counting multiplicity) which satisfy $p\equiv a\bmod q$.
Also, biases exhibited in families of weighted sums, such as $L_\alpha(x)=\sum_{n=1}^{\floor{x}} \lambda(n)/n^\alpha$  and $M_\alpha(x)=\sum_{n=1}^{\floor{x}} \mu(n)/n^\alpha$ for varying real $\alpha$, were studied in \cite{Alkan20,Alkan21,MosTru,MT2017}.

\section{Functions with values in $\{-1,1\}$ and the proof of Theorem~\ref{thmpm1}}\label{secpm1}

Suppose $\fake(n)$ is a multiplicative function with $\fake(p^k)=\epsilon_k\in\{-1,1\}$ for every prime $p$ and each $k\geq1$, with $\epsilon_k$ independent of $p$.
Clearly $\epsilon_0=1$. For $\sigma>1$ write
\[
T(s) = \sum_{n\geq1} \frac{\fake(n)}{n^s}.
\]

Suppose first that $\epsilon_1=1$.
Then we can approximate $T(s)$ with the Riemann zeta function:
\[
T(s) = \prod_p\left(1+\frac{1}{p^s}+\cdots\right) \approx \prod_p\left(1-\frac{1}{p^s}\right)^{-1} = \zeta(s).
\]
We thus write $T(s)=\zeta(s)U(s)$ and solve for $U(s)$.
If $\epsilon_2=1$, then
\[
U(s) = \prod_p \left(1+\frac{1}{p^s}+\frac{1}{p^{2s}}+\frac{\epsilon_3}{p^{3s}}+\cdots\right)\left(1-\frac{1}{p^s}\right) = \prod_p\left(1+ \frac{\epsilon_3-1}{p^{3s}} + \cdots\right),
\]
so $U(s)$ is analytic on $\sigma>1/3$, and $T(s)$ has no pole at $s=1/2$.
In this case $T(s)$ does have a pole at $s=1$, so the expression for $\Fake(x)/\sqrt{x}$ from Perron's formula will exhibit a linear term $ax$ (with $a=U(1)$).
However this has no effect after our normalization \eqref{eqnNormalize}, so such $\Fake(x)$ exhibit no persistent or apparent bias.

If $\epsilon_2=-1$, then
\begin{align*}
U(s) &= \prod_p \left(1+\frac{1}{p^s}-\frac{1}{p^{2s}}+\cdots\right)\left(1-\frac{1}{p^s}\right) = \prod_p\left(1- \frac{2}{p^{2s}} + \cdots\right)\\
&\approx \prod_p\left(1-\frac{1}{p^s}\right)^2 = \frac{1}{\zeta(2s)^2}.
\end{align*}
Let $U_1(s) = T(s) \zeta(2s)^2/\zeta(s)$, so
\[
U_1(s) = \prod_p \frac{1-p^{-s}}{\left(1-p^{-2s}\right)^2}\left(1+\frac{1}{p^s}-\frac{1}{p^{2s}}+\frac{\epsilon_3}{p^{3s}}+\cdots\right) = \prod_p\left(1+ \frac{1+\epsilon_3}{p^{3s}} + \cdots\right),
\]
and so $U_1(s)$ converges for $\sigma>1/3$.
While $T(s) = U_1(s)\zeta(s)/\zeta(2s)^2$ in general has a pole at $s=1$ with residue $U_1(1)/\zeta(2)^2$, it is analytic at $s=1/2$, so again $\Fake(x)$ exhibits no persistent or apparent bias.

Suppose next that $\epsilon_1=\epsilon_2=-1$, so this includes the case where $\fake(n)=\xi(n)$.
Then we approximate $T(s)$ using $1/\zeta(s)$, and write $U(s)=T(s)\zeta(s)$, so
\[
U(s) = \prod_p \left(1-\frac{1}{p^{2s}}+\cdots\right) \approx \frac{1}{\zeta(2s)}.
\]
We thus set $U_1(s)=U(s)\zeta(2s)$, so
\[
U_1(s) = T(s)\zeta(s)\zeta(2s) = \prod_p \left(1 + \frac{\epsilon_3-1}{p^{3s}} + \cdots\right),
\]
which converges on $\sigma>1/3$.
Thus $T(s) = U_1(s)/\zeta(s)\zeta(2s)$ has no pole at $s=1/2$ (and furthermore $a=0$ here as well), so $\Fake(x)$ exhibits no bias $b$.

We turn then to the remaining case, where $\epsilon_1=-1$ and $\epsilon_2=1$, which includes the case $\fake(n)=\lambda(n)$.
We have
\[
T(s) = \prod_p\left(1-\frac{1}{p^s}+\frac{1}{p^{2s}}+\cdots\right) \approx \prod_p\left(1+\frac{1}{p^2}\right)^{-1} = \frac{\zeta(2s)}{\zeta(s)},
\]
so we write $U(s) = T(s)\zeta(s)/\zeta(2s)$, and calculate
\[
U(s) = \prod_p\left(1 + \frac{1+\epsilon_3}{p^{3s}} + \cdots\right),
\]
so $U(s)$ converges on $\sigma>1/3$.
Thus $T(s)=U(s)\zeta(2s)/\zeta(s)$ has a pole at $s=1/2$ with residue $U(1/2)/2\zeta(1/2)$, provided $U(1/2)\neq0$, so from Perron's formula the constant term for $\Fake(x)/\sqrt{x}$ is $U(1/2)/\zeta(1/2)$, and this accounts for the apparent bias of $\Fake(x)$.
We still have freedom in the selection of $\epsilon_k$ for $k\geq3$, and we can choose these to manipulate the bias. We have
\begin{equation}\label{eqnU4}
U(s) = \prod_p\left(1 + \sum_{k\geq3} \frac{\epsilon_{k-1}+\epsilon_k}{p^{sk}}\right) =: \prod_p C_p(s).
\end{equation}
For each prime $p$, we have
\[
1 - \sum_{k\geq4} \frac{2}{p^{k/2}} \leq C_p(1/2) \leq 1 + \sum_{k\geq3} \frac{2}{p^{k/2}},
\]
with the lower bound achieved by selecting $\epsilon_k=-1$ for $k\geq3$, and the upper bound by choosing $\epsilon_k=1$ for $k\geq3$.
A straightforward calculation then verifies that the bias of $\Fake(x)$ lies in the interval $[A_1,B_1]$ of Theorem~\ref{thmpm1}.
We note that the minimal value $A_1=-10.29174\ldots$ is attained when $\epsilon_k=-1$ if and only if $k=1$ as in \eqref{eqnfakemin}, and the maximal value $B_1=0.16918\ldots$ is achieved when $\epsilon_k=1$ if and only if $k\in\{0,2\}$.

It remains to show that there exist uncountably many $\fake(n)$ with no apparent bias having $\epsilon_1=-1$ and $\epsilon_2=1$, and to characterize these functions.
For $p\geq3$, we have
\[
C_p(1/2) \geq 1 - \frac{2}{p^2(1-p^{-1/2})} \geq 1 - \frac{2}{9(1-3^{-1/2})} > 0.47,
\]
so $U(1/2)\geq0$ (with $U(s)$ in \eqref{eqnU4}) if and only if $C_2(1/2)\leq0$.
Each $C_p(s)$ converges for $\sigma>0$, so to determine the sign of $C_2(1/2)$ it is enough to determine the sign of the $p=2$ factor in the Euler product for $T(s)$ at $s=1/2$.
Denote this value by $\gamma_2$.
For $k\geq3$ set $\delta_k\in\{0,1\}$ so that $\epsilon_k=1-2\delta_k$.
Then we calculate
\[
\gamma_2 = 1 - \frac{1}{2^{1/2}} + \frac{1}{2} + \sum_{k\geq3} \frac{1-2\delta_k}{2^{k/2}}
= 2\left(1-\sum_{k\geq3} \frac{\delta_k}{2^{k/2}}\right).
\]
Define two real numbers $\alpha$ and $\beta$ in $[0,1]$ by writing their binary expansions using the values of $\delta_k$ in the following way:
\begin{equation}\label{eqnAlphaBeta}
\alpha = (0.\delta_3 \delta_5 \delta_7\ldots)_2,\quad
\beta = (0.\delta_4 \delta_6 \delta_8\ldots)_2.
\end{equation}
Then
\[
\gamma_2 = 2\left(1-\frac{\alpha}{\sqrt{2}} - \frac{\beta}{2}\right),
\]
and $U(1/2)\geq0$ precisely when $\alpha\sqrt{2}+\beta \geq 2$.
Thus we have zero apparent bias in this case precisely when
$\alpha\sqrt{2}+\beta = 2.$
There are uncountably many ways to choose the $\delta_k$ to ensure this, arising from the real solutions to $\alpha\sqrt{2}+\beta = 2$ with $\alpha$ and $\beta$ in $[0,1]$.
Note that if $\alpha$ or $\beta$ is a dyadic rational in $(0,1)$, then it has two relevant binary expansions, and either one may be used to construct a function with no apparent bias when $\alpha\sqrt{2}+\beta=2$. 
This completes the proof of Theorem~\ref{thmpm1}.\qed

\vskip\baselineskip

\begin{figure}[tb]
\begin{center}
\includegraphics[width=4in]{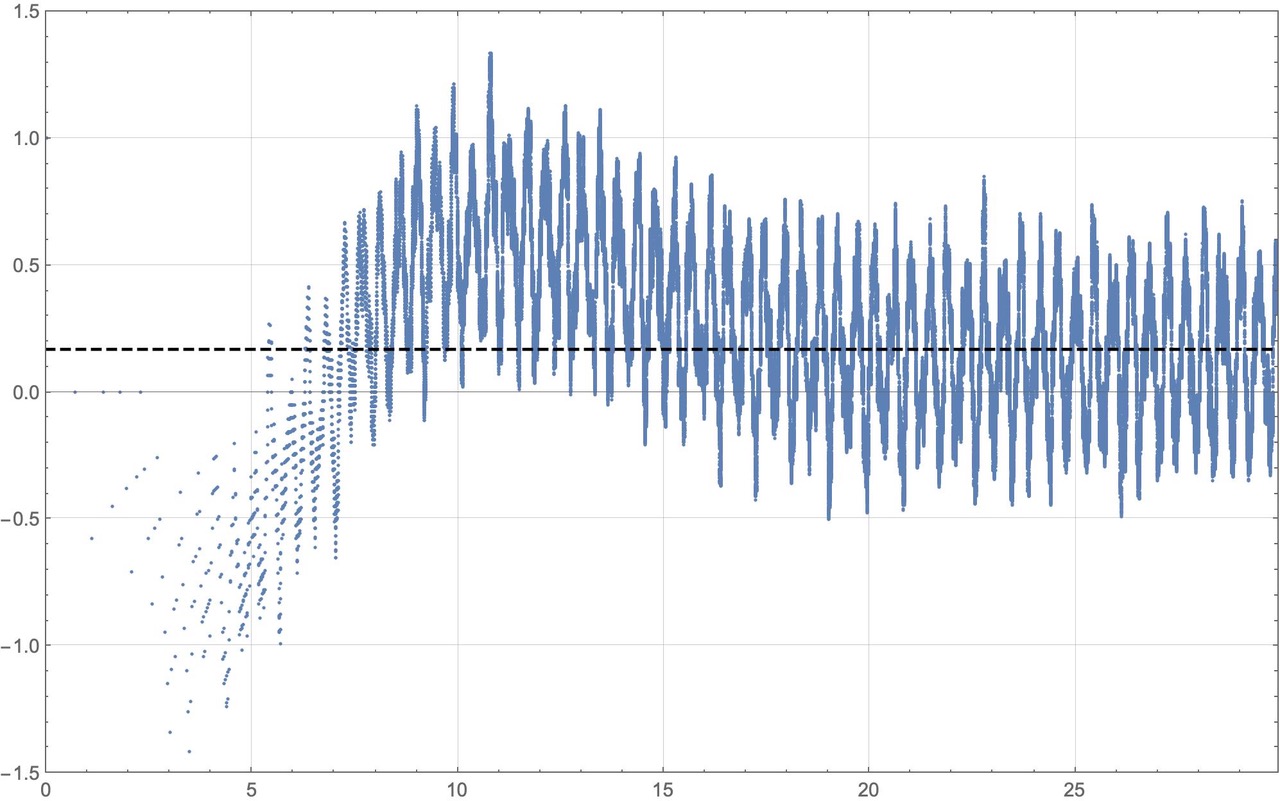}
\end{center}
\caption{Normalized summatory function $e^{-u/2}\Fake_{\textrm{max}}(e^u)$ with $u\leq\log(10^{13})$, for $\fake_{\textrm{max}}(n)$ from \eqref{eqnfakemax}, with maximal apparent bias.}\label{figMaxPosFakeLambdas}
\end{figure}

We can describe the functions achieving the extremal values in part~(\ref{item1b}) of Theorem~\ref{thmpm1}.
For a positive integer $n$, write $n = \prod_k n_{(k)}$, where
\begin{equation}\label{eqnnk}
n_{(k)} = \prod_{p^k\parallel n} p^k.
\end{equation}
Then the minimal apparent bias in this family is attained by the function $\fake_{\textrm{min}}(n)=\mu(n_{(1)})=\lambda(n_{(1)})=\xi(n_{(1)})$ as in \eqref{eqnfakemin}, and the maximal one by
\begin{equation}\label{eqnfakemax}
\fake_{\textrm{max}}(n)=\xi(n/n_{(2)}).
\end{equation}
Plots showing the values for the normalized summatory functions for $\fake_{\textrm{min}}(n)$ and $\fake_{\textrm{max}}(n)$ up to $10^{13}$ are shown respectively in Figures~\ref{figMaxNegBias} and~\ref{figMaxPosFakeLambdas}.

\begin{figure}[tb]
\begin{center}
\includegraphics[width=4in]{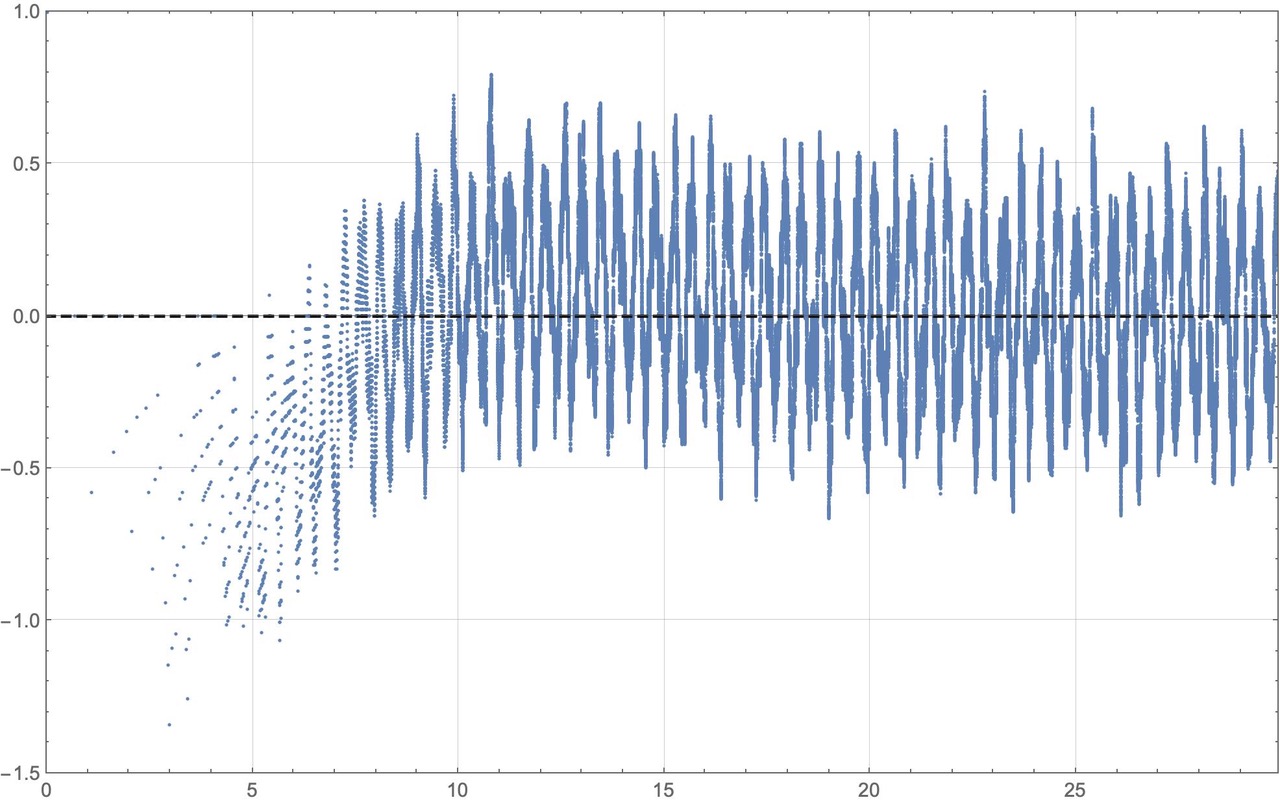}
\end{center}
\caption{Normalized summatory function $e^{-u/2}\Fake_{(1/\sqrt{2},1)}(e^u)$ with $u\leq\log(10^{13})$, for $\fake_{(1/\sqrt{2},1)}(n)$ from the end of Section~\ref{secpm1}, with $\epsilon_1=-1$, $\epsilon_2=1$, and no apparent bias.}\label{figNoBias}
\end{figure}

We can also construct an example of a function from part~(\ref{item1c}) of Theorem~\ref{thmpm1} with no apparent bias, by selecting $\alpha=1/\sqrt{2}$ and $\beta=1$ in \eqref{eqnAlphaBeta}.
That is, we select the values $\delta_{2j+1}$ for $j\geq1$ according to the binary expansion of $1/\sqrt{2} = (0.1011010100\ldots)_2$, and set the values $\delta_{2j}=1$ for $j\geq2$.
We denote this function by $\fake_{(1/\sqrt{2},1)}(n)$.
A plot showing values up to $10^{13}$ for its normalized summatory function appears in Figure~\ref{figNoBias}.

\section{Functions with values in $\{-1,0,1\}$ and the proof of Theorem~\ref{thmzpm1}}\label{seczpm1}

The proof of Theorem~\ref{thmzpm1} proceeds in the same fashion.
The arguments covering the unbiased cases of Theorem~\ref{thmpm1} lift directly to the current environment, where $\epsilon_k=0$ is allowed as well, so we check only the new cases claimed here.
First, if $\epsilon_1=1$ and $\epsilon_2=0$, then
\[
T(s) = \prod_p \left(1+\frac{1}{p^s}+\frac{\epsilon_3}{p^{3s}}+\cdots\right) \approx \frac{\zeta(s)}{\zeta(2s)},
\]
so we define $U(s)=T(s)\zeta(2s)/\zeta(s)$, and we check that $U(s)$ converges for $\sigma>1/3$.
Thus $T(s)=U(s)\zeta(s)/\zeta(2s)$ has no pole at $s=1/2$, and the corresponding summatory function $\Fake(x)$ is not biased at the $\sqrt{x}$ scale.
(There is a potential pole at $s=1$, so $a = U(1)/\zeta(2)$ in this case.)

Next, suppose $\epsilon_1=\epsilon_2=0$.
Clearly we are in the unbiased case if $\epsilon_k=0$ for all $k\geq3$, so we may assume some $\epsilon_k\neq0$.
Let $k_0\geq3$ be the smallest positive such index.
Then
\[
T(s) = \prod_p \left(1+\frac{\epsilon_{k_0}}{p^{k_0 s}}+\cdots\right) \approx \zeta(k_0 s)^{\epsilon_{k_0}}.
\]
Set $U(s)=T(s)/\zeta(k_0 s)^{\epsilon_{k_0}}$, so that $U(s)$ converges on $\sigma>1/(k_0+1)$, so $\Fake(x)$ has no persistent or apparent bias.
Similarly, if $\epsilon_1=0$ and $\epsilon_2=-1$, then
\[
T(s) = \prod_p \left(1-\frac{1}{p^{2s}}+\cdots\right) \approx \frac{1}{\zeta(2s)},
\]
and $U(s)=T(s)\zeta(2s)$ converges when $\sigma>1/3$, so $T(s)$ is analytic at $s=1/2$ (in fact, $T(1/2) = 0$), and the corresponding summatory function has no bias of either type.
Last, assume $\epsilon_1=-1$ and $\epsilon_2=0$, which includes the case  $\fake(n)=\mu(n)$.
Then
\[
T(s) = \prod_p \left(1-\frac{1}{p^s}+\frac{\epsilon_3}{p^{3s}}+\cdots\right) \approx \frac{1}{\zeta(s)},
\]
and one may check that $U(s)=T(s)\zeta(s)$ converges for $\sigma>1/3$, and $T(s)$ has no pole at $s=1/2$.

If $\epsilon_1=-1$ and $\epsilon_2=1$, then much of the analysis in Section~\ref{secpm1} for this case generalizes easily, and in particular the extremal solutions are not disturbed by the possibility that some $\epsilon_k=0$.
The functions covered by this case include the $\mu_r(n)$ with $r\geq3$ studied by Tanaka.
We can characterize the functions in this family having no apparent bias in a similar way, by determining when the factor with $p=2$ in the Euler product for $T(s)$ is $0$ at $s=1/2$.
For $k\geq3$, write $\epsilon_k=\delta_k^+-\delta_k^-$, with $\delta_k^+,\delta_k^-\in\{0,1\}$ and $\delta_k^+\delta_k^-=0$.
Then no apparent bias occurs precisely when
\[
1 - \frac{1}{\sqrt{2}} + \frac{1}{2} + \sum_{k\geq3} \frac{\delta_k^+-\delta_k^-}{2^{k/2}} = 0,
\]
that is, when
\[
\frac{3}{2} - \frac{1}{\sqrt{2}} + \frac{1}{\sqrt{2}}\sum_{k\geq1} \frac{\delta_{2k+1}^+-\delta_{2k+1}^-}{2^k} + \frac{1}{2}\sum_{k\geq1} \frac{\delta_{2k+2}^+-\delta_{2k+2}^-}{2^k} = 0.
\]
Let $\alpha^+$, $\alpha^-$, $\beta^+$, and $\beta^-$ be real numbers in $[0,1]$ defined by
\begin{align*}
&\alpha^+ = (0.\delta_3^+\delta_5^+\delta_7^+\ldots)_2, \quad
\beta^+ = (0.\delta_4^+\delta_6^+\delta_8^+\ldots)_2, \\
&\alpha^- = (0.\delta_3^-\delta_5^-\delta_7^-\ldots)_2, \quad
\beta^- = (0.\delta_4^-\delta_6^-\delta_8^-\ldots)_2,
\end{align*}
and let $\alpha=\alpha^+-\alpha^-$ and $\beta=\beta^+-\beta^-$.
Then no apparent bias in this family occurs precisely when
$2\alpha + \sqrt{2}\beta = 2 - 3\sqrt{2}.$

We turn then to the remaining case, where $\epsilon_1=0$ and $\epsilon_2=1$.
Then
\[
T(s) = \prod_p \left(1+\frac{1}{p^{2s}}+\frac{\epsilon_3}{p^{3s}}+\cdots\right) \approx \zeta(2s),
\]
so we set $U(s)=T(s)/\zeta(2s)$ and verify that it converges for $\sigma>1/3$.
We see that $T(s)=U(s)\zeta(2s)$ has a pole at $s=1/2$ with residue $U(1/2)/2$ (provided $U(1/2)\neq0$), so the constant term in the expression for $\Fake(x)/\sqrt{x}$ from Perron's formula is $U(1/2)$.
Further, since $T(s)$ is analytic in $\sigma>1/3$ except for the simple pole at $s=1/2$, it follows that $\Fake(x)/\sqrt{x}\to U(1/2)$ as $x\to\infty$, so $\Fake(x)$ has a persistent bias.
We can select $\epsilon_k$ with $k\geq3$ to manipulate its value.
We compute
\begin{align*}
U(1/2) &= \prod_p \left(1+\sum_{k\geq3} \frac{\epsilon_k-\epsilon_{k-2}}{p^{k/2}}\right)\\
&= \prod_p \left(1+\left(1-\frac{1}{p}\right)\left(\sum_{k\geq3} \frac{\epsilon_k}{p^{k/2}}\right) - \frac{1}{p^2}\right) =: \prod_p Q_p.
\end{align*}
Now for each prime $p$ we have
\[
Q_p \geq 1+\left(1-\frac{1}{p}\right)\left(\sum_{k\geq3} \frac{-1}{p^{k/2}}\right) - \frac{1}{p^2} = 1-\frac{1}{p^{3/2}}-\frac{2}{p^2} \geq \frac{1}{2}-\frac{1}{2\sqrt{2}} > 0.14,
\]
and
\[
Q_p \leq 1+\left(1-\frac{1}{p}\right)\left(\sum_{k\geq3} \frac{1}{p^{k/2}}\right) - \frac{1}{p^2} = 1+\frac{1}{p^{3/2}},
\]
so
\[
\prod_p \left(1-\frac{1}{p^{3/2}}-\frac{2}{p^2}\right) \leq U(1/2) \leq \frac{\zeta(3/2)}{\zeta(3)}.
\]
This completes the proof of Theorem~\ref{thmzpm1}.\qed

\vskip\baselineskip

\begin{figure}[tb]
\begin{center}
\includegraphics[width=4in]{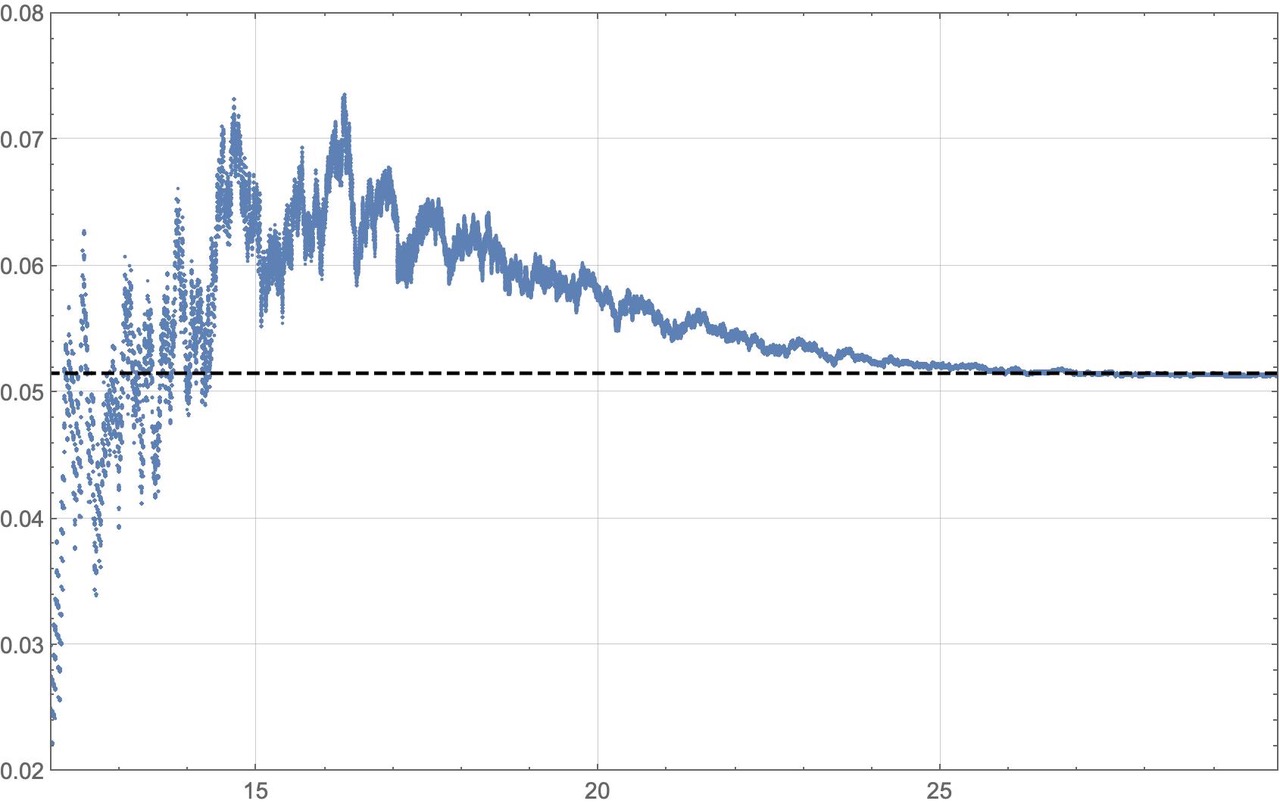}
\end{center}
\caption{Normalized summatory function $e^{-u/2}\Fake_{\textrm{Min}}(e^u)$ with $u\leq\log(10^{13})$, for $\fake_{\textrm{Min}}(n)$ from \eqref{eqnfakeMin}, with minimal persistent bias.}\label{figMinBias01}
\end{figure}

We remark that the maximal persistent bias in the proof of Theorem~\ref{thmzpm1} occurs for the function $\fake_{\textrm{Max}}(n)$ from \eqref{eqnfakeMax}, which is $1$ when $n$ has no prime factors with multiplicity $1$ and $0$ otherwise, that is, when $\fake_{\textrm{Max}}(n)$ is the indicator function for the \textit{powerful} integers.
The minimal persistent bias also occurs when $\fake(n)$ is supported on the powerful integers: using \eqref{eqnnk} the extremal function is
\begin{equation}\label{eqnfakeMin}
\fake_{\textrm{Min}}(n) = \begin{cases}
\xi(n/n_{(2)}), & n_{(1)} = 1,\\
0, & n_{(1)}>1.
\end{cases}
\end{equation}
The functions $\Fake_{\textrm{Max}}(x)/\sqrt{x}$ and $\Fake_{\textrm{Min}}(x)/\sqrt{x}$ for $x\leq10^{13}$ are displayed in Figures~\ref{figMaxPosBias} and~\ref{figMinBias01} respectively.
Note that the Dirichlet series for $\fake_{\textrm{Max}}(n)$ is $T_{\textrm{Max}}(s) = \zeta(2s)\zeta(3s)/\zeta(6s)$, which (on RH) has simple poles on the $\sigma=1/12$ line, in addition to the ones at $s=1/2$ and $s=1/3$, so from Perron's formula we expect the oscillations in the normalized function $\Fake_{\textrm{Max}}(x)/\sqrt{x}$ that arise from the complex poles to attenuate like $x^{-5/12}$ as $x$ grows large.
Also, one may verify that the Dirichlet series for $\fake_{\textrm{Min}}(n)$ may be written as $U_1(s)\zeta(2s)/\zeta(3s)\zeta(4s)^2$, with $U_1(s)$ analytic on $\sigma>1/7$, so (again assuming RH) the oscillations in $\Fake_{\textrm{Min}}(x)/\sqrt{x}$ arising from the simple poles of $T_{\textrm{Min}}(s)$ on the $\sigma=1/6$ line will dampen as $x^{-1/3}$ as $x\to\infty$.

\section{Proof of Theorem~\ref{thmBias}}\label{secBias}

Let $\fake(n)$ be an arithmetic function with the properties in the hypothesis, and let $T(s)$ denote its Dirichlet series.
Then $T(s)=U(s)\zeta(2s)/\zeta(s)$ with $U(s)$ given by \eqref{eqnU4}, and the apparent bias of $\Fake(x)$ is $b=U(1/2)/\zeta(1/2)$.
Let $r$ be a real constant whose value will be selected later.
It follows from Abel's summation formula that
\[
\zeta(2s) = s\int_1^\infty \frac{\floor{\sqrt{x}}}{x^{s+1}}\,dx
\]
and
\begin{equation}\label{eqnFakeIntegral}
T(s) - (b+r)\zeta(2s) = s\int_1^\infty \frac{\Fake(x) - (b+r)\floor{\sqrt{x}}}{x^{s+1}}\,dx
\end{equation}
for $\sigma>1$.
Suppose that this integrand has constant sign for all sufficiently large $x$.
Then it follows from Landau's theorem (see \cite[Lem.\ 15.1]{MV}) that this function is analytic in $\sigma>1/2$, and thus so is $T(s)$, and RH follows.
Let $\rho_1 = 1/2+i\gamma_1$ denote the first zero of the zeta function on the critical line.
Using \eqref{eqnFakeIntegral}, we have
\[
\abs{T(\sigma+i\gamma_1) - (b+r)\zeta(2\sigma+2i\gamma_1)} < 2\abs{\rho_1}\abs{T(\sigma) - (b+r)\zeta(2\sigma)}
\]
for $\sigma>1/2$, so
\begin{equation}\label{eqnTIneq}
\begin{split}
&\frac{1}{\abs{\rho_1}}\lim_{\sigma\to\frac{1}{2}^+} \left(\sigma-\frac{1}{2}\right)\abs{T(\sigma+i\gamma_1) - (b+r)\zeta(2\sigma+2i\gamma_1)}\\
&\qquad \leq
2\lim_{\sigma\to\frac{1}{2}^+} \left(\sigma-\frac{1}{2}\right) \abs{T(\sigma) - (b+r)\zeta(2\sigma)}.
\end{split}
\end{equation}
Since
\[
\lim_{\sigma\to\frac{1}{2}^+} \left(\sigma-\frac{1}{2}\right) \zeta(2\sigma) = \frac{1}{2}
\quad\textrm{and}\quad
\lim_{\sigma\to\frac{1}{2}^+} \abs{T(\sigma) - b\zeta(2\sigma)} < \infty,
\]
it follows that
\[
\lim_{\sigma\to\frac{1}{2}^+} \left(\sigma-\frac{1}{2}\right) \bigl(T(\sigma) - (b+r)\zeta(2\sigma)\bigr) = \frac{r}{2}.
\]
In addition,
\begin{equation}\label{eqnLim}
\lim_{\sigma\to\frac{1}{2}^+} \left(\sigma-\frac{1}{2}\right) \bigl(T(\sigma+i\gamma_1) - (b+r)
\zeta(2\sigma+2i\gamma_1)\bigr) = \lim_{\sigma\to\frac{1}{2}^+} \left(\sigma-\frac{1}{2}\right) T(\sigma+i\gamma_1).
\end{equation}
Since $\rho_1$ is a simple zero of $\zeta(s)$, the limit in \eqref{eqnLim} is finite.
Denote its value by $c$.
Given that
\[
\frac{1}{U(s)} = \prod_p \left(1-\frac{1+\epsilon_3}{p^{3s}}+\cdots\right),
\]
we have that $1/U(s)$ is analytic on $\sigma>1/3$, so $U(\rho_1)\neq0$ and consequently $c\neq0$.
Then from \eqref{eqnTIneq} we have
\[
\abs{\frac{c}{\rho_1}} \leq \abs{r} .
\]
By selecting $r = \abs{c/2\rho_1}$ so that $\abs{r}<\abs{c/\rho_1}$, we conclude that each of $\Fake(x) - (b\pm r)\floor{\sqrt{x}}$ must change sign infinitely often as $x\to\infty$.
The theorem follows.
\qed

\vskip\baselineskip

Finally, we remark that Ingham's observation regarding $M(x)$ and $L(x)$ may be modified to show that if $\Fake_{\textrm{min}}(x)/\sqrt{x}$ as in Figure~\ref{figMaxNegBias} is bounded either above or below by some constant, then the Riemann hypothesis follows, as well as the simplicity of the zeros of the Riemann zeta function.
This follows by an argument similar to that employed in the prior proof.
By analogy with P\'olya's question for the sum of the Liouville function, can one show that $\Fake_{\textrm{min}}(x)$ has infinitely many sign changes?
Can one determine an $x\geq42$ where $\Fake_{\textrm{min}}(x)>0$?
Owing to the computations used to make Figure~\ref{figMaxNegBias}, any such $x$ must satisfy $x>10^{13}$.

\section*{Acknowledgements}

This research was undertaken with the assistance of resources and services from the National Computational Infrastructure (NCI), which is supported by the Australian Government.

 \end{document}